\documentclass[11pt]{article}
\usepackage{graphicx,psfrag}
\usepackage{amssymb,amsfonts,amsthm}
\usepackage{amsmath,amsthm,amssymb}
\usepackage{amsfonts}

\newcommand{\HNN}{{\mathrm HNN}}
\newcommand{\SL}{{\mathrm SL}}
\newcommand{\la}{\langle}
\newcommand{\ra}{\rangle}

\newtheorem{example}{Example}[section]

\newtheorem{theorem}{Theorem}[section]
\newtheorem{problem}[theorem]{Problem}

\newtheorem{definition}[theorem]{Definition}

\newtheorem{prob}[theorem]{Problem}
\newcommand{\lo }{\lim^\omega }
\newcommand{\oom }{o_\omega }
\newcommand{\Oo }{O_\omega }
\newcommand{\To }{\Theta _\omega }

\newcommand{\ch}{lacunary hyperbolic\,}

\newcommand{\CG }{{\rm Con} ^{\omega } (G, d)}

\newcommand{\Con}{{\mathrm{Con}}}
\newcommand{\lm}{{\lim}}
\newcommand{\dist}{{\mathrm{dist}}}

\newcommand{\pp}{{\mathcal P}}

\newcommand {\N}{\mathbb{N}} %% positive integers
\newcommand {\Z}{\mathbb{Z}}            %% integers
\newcommand {\R}{\mathbb{R}} %% reals
\newcommand {\iv}{^{-1}}
\newcommand{\lio}[1]{\lm^\omega\left(#1\right)}

\begin{document}

\title{Some group theory problems}
\author{Mark Sapir\thanks{The author was supported in part by the NSF grants
DMS 0245600, DMS 0455881 and a BSF (the US-Israeli) grant.}}
\date{}
\maketitle

\qquad\qquad\qquad\qquad\qquad\qquad\qquad\qquad\qquad\qquad
Dedicated to Boris Isaakovich Plotkin

\begin{abstract} This is a survey of some problems in geometric
group theory which I find
interesting. The problems are from different areas of group theory.
Each section is devoted to problems in one area. It contains an
introduction where I give some necessary definitions and
motivations, problems and some discussions of them. For each
problem, I try to mention the author. If the author is not given,
the problem, to the best of my knowledge, was formulated by me
first.
\end{abstract}

\tableofcontents

\section{Groups with extreme properties (monsters)}
\label{s1}

I assume that the reader knows basic definitions of group theory, as
well as the definitions of amenable group, and the growth function
of a group.

 There are currently several methods of constructing
groups with unusual properties (say, Tarski monsters, i.e. groups
with all proper subgroups cyclic, or non-Amenable groups without
free non-Abelian subgroups, or finitely generated infinite groups
with finitely many conjugacy classes). The methods were originated
in works by Adian, Novikov, Olshanskii \cite{Ad,book}, Gromov
\cite{Gr,Gr3} , and developed by their students. The latest advance
in the ``theory of monsters" was due to Osin who constructed an
infinite finitely generated group with exactly two conjugacy classes
\cite{OsinTwoClasses}. Each of these constructions proceeds as
follows: we start with a free non-Abelian group (or a non-Abelian
hyperbolic group in general), and then add relations one by one, at
each step solving a little portion of the problem. It is important
that at each step, we get a hyperbolic (or relatively hyperbolic)
group, so that we can proceed by induction. In the limit, we get a
group satisfying the desired property. For example, if we construct
an infinite finitely generated group satisfying the law $x^n=1$,
then (following the method from \cite{book}), at each step we add
one relation of the form $u^n=1$, so that in the limit the law
$x^n=1$ is satisfied for all $x$.

The groups constructed this way are usually infinitely presented:
relations at each step are chosen so that they ``do not interfere"
with the relations from the previous step (we usually need some kind
of small cancelation condition to achieve that). Hence relations
added later do not follow from the previous relations. In this
section, by ``monsters", I understand finitely generated infinite
groups satisfying any of the following properties. These are only
samples of possible properties, one can easily find generalizations
or similar properties which are also of great interest.

\begin{itemize}
\item[(1)] The group is non-amenable but contains no non-Abelian
free subgroups (solving the so-called Day-von Neumann problem
\cite{book});
\item[(2)] The group satisfies the law $x^n=1$ (Burnside problem, \cite{Ad}, \cite{book});
\item[(3)] Every proper subgroup is cyclic (Tarski problem \cite{book});
\item[(4)] Every proper subgroup is finite and cyclic of the same fixed prime order
(Tarski problem \cite{book});
\item[(5)] The group is complete (i.e. every element has roots of
any degree, Guba, \cite{Gu}, Mikhajlovskii, Olshanskii \cite{MO});
\item[(6)] The number of conjugacy classes is finite (S. Ivanov,
\cite{book}), preferably 2 \cite{OsinTwoClasses}.
\end{itemize}
I would add to this list the following property which so far is
dealt with by completely different methods (Grigorchuk,
\cite{Grigorchuk}, \cite{GNS}), but all known examples are also
infinitely presented. The property is very important:

\begin{itemize}
\item[(7)] The group has intermediate -- not bounded above by a
polynomial and not bounded below by an exponent -- growth function
(Milnor problem).
\end{itemize}

The following problem is very old, and probably was first formulated
by Kurosh for torsion groups in the 50s, or even before that by
somebody else. I heard the problem for other types of monsters from
Adian, Grigorchuk, Olshanskii and others.

\begin{problem} Find finitely presented monsters of types (1)-(7) or prove that they
do not exist.
\end{problem}

So far only finitely presented monsters of type (1) (non-Amenable
but without non-Abelian free subgroups) have been found
\cite{OSamen}. There are some ideas how to construct monsters of
type (2). In \cite{OSamen}, Olshanskii and I constructed a finitely
presented bounded torsion group $G$ with finite endomorphic
presentation (S. Ivanov later constructed another example in
\cite{Ivanov}): the group is generated by a finite set $S$ and there
exists an injective endomorphism $\phi$ of the free group with basis
$S$ and a finite list of relations $r_1,...,r_n$ of the group such
that the words $\phi^m(r_i)$ for a presentation of $G$ ($m=0,1,...;
i=1,...,n$). Finite endomorphic presentations of groups are
interesting because such a group $G$ has a finitely presented
ascending HNN extension $\HNN_\phi(G)$, and ascending HNN extensions
preserve many properties of the group including amenability and
elementary amenability. Ascending HNN extensions and endomorphic
presentations have been used to produce finitely presented solvable
groups by Remeslennikov \cite{Rem} and Baumslag \cite{Bau}. Lysenok
discovered such a presentation for Grigorchuk's group. Many other
groups acting on locally finite rooted trees have such presentations
(see the survey Grigorchuk-Nekrashevich-Sushchanskii \cite{GNS}).

Although I am reasonably sure that a finitely presented infinite
bounded torsion group exists (note that unbounded torsion infinite
finitely presented groups are not known also), I am not at all sure
that finitely presented monsters of types (3)-(7) exist. There are
results which show that finitely presented groups satisfy nicer
properties than arbitrary finitely generated groups. One example of
such a result is the theorem of M. Kapovich and Kleiner: if one
asymptotic cone of a finitely presented group is a tree, then the
group is hyperbolic (see the appendix of \cite{OOS}). In \cite{OOS},
Olshanskii, Osin and I showed that for finitely generated groups
this is is extremely far from being true.

It is interesting that (as I have mentioned above) Grigorchuk's
group has also a finite endomorphic presentation, and is torsion.
But it does not have bounded torsion. It is very interesting to know
whether there are infinite bounded torsion groups of intermediate
growth, or, more generally, whether there are finitely generated
infinite amenable bounded torsion groups. Note that all torsion
automata groups (groups constructed by Sushchanskii, Grigorchuk,
Gupta-Sidki and others \cite{GNS}) are residually finite by
construction (they act faithfully on rooted locally finite trees),
but there are no bounded torsion residually finite groups. That
follows from Zelmanov's solution of the restricted Burnside problem
\cite{Ze1,Ze2}.

The following even more general problem is again a part of old
folklore.

\begin{problem} Are there (finitely generated non-virtually cyclic)
amenable monsters of types (2)-(6).
\end{problem}

Note that groups of intermediate growth are all amenable by the
F\"olner criterium (F\"olner sets are simply balls in that case).

\section{Asymptotic cones of groups}

Asymptotic cones were introduced by Gromov in \cite{Gro81}, a
definition via ultrafilters was given by van den Dries and Wilkie
\cite{DW}. An asymptotic cone of a metric space is, roughly
speaking, what one sees when one looks at the space from infinitely
far away.

Here is a more precise definition of asymptotic cones. A
non-principal ultrafilter $\omega$ is a finitely additive measure
defined on all subsets $S$ of ${\mathbb N}$, such that $\omega(S)
\in \{0,1\}$, $\omega (\mathbb N)=1$, and $\omega(S)=0$ if $S$ is a
finite subset. An ultrafilter allows one to define the limit of any
bounded sequence of numbers $x_n$, $n\in \N$: the limit $\lo x_n$
with respect to $\omega$ is the unique real number $a$ such that
$\omega(\{i\in {\mathbb N}: |x_i-a|<\epsilon\})=1$ for every
$\epsilon>0$. Similarly, $\lo x_n=\infty $ if $\omega(\{i\in
{\mathbb N}: |x_i|>M \})=1$ for every $M>0$.

Given two infinite sequences of real numbers $(a_n)$ and $(b_n)$ we
write $a_n=\oom (b_n)$ if $\lo a_n/b_n =0$. Similarly $a_n=\To
(b_n)$ (respectively $a_n=\Oo (b_n)$) means that $0<\lo
(a_n/b_n)<\infty $ (respectively $\lo (a_n/b_n)<\infty $).

Let $(X_n,\dist_n)$, $n\in\N$, be a metric space. Fix an arbitrary
sequence $e=(e_n)$ of points $e_n\in X_n$. Consider the set ${\cal
F}$ of sequences $g=(g_n)$, $g_n\in X_n$, such that
$\dist_n(g_n,e_n) \le c $ for some constant $c=c(g)$. Two sequences
$(f_n)$ and $(g_n)$ of this set ${\cal F}$ are said to be {\em
equivalent} if $\lio{\dist_n(f_n,g_n)} =0$. The equivalence class of
$(g_n)$ is denoted by $(g_n)^\omega$. The $\omega$-{\em limit}
$\lio{ (X_n)_e}$ is the quotient space of equivalence classes where
the distance between $(f_n)^\omega$ and $(g_n)^\omega$ is defined as
$\lio{\dist(f_n,g_n)}$.

An {\em asymptotic cone} $\Con^\omega(X,e, d)$ of a metric space
$(X,\dist)$ where $e=(e_n)$, $e_n\in X$, and $d=(d_n)$ is an
unbounded non-decreasing {\it scaling sequence} of positive real
numbers, is the $\omega$-limit of spaces $X_n=(X,\dist/d_n)$. The
asymptotic cone is a complete space; it is a geodesic metric space
if $X$ is a geodesic metric space (\cite{Gr3,Drutu}). Note that
$\Con^\omega(X, e,d)$ does not depend on the choice of $e$ if $X$ is
homogeneous (say, if $X$ is a finitely generated group with a word
metric), so in that case, we shall omit $e$ in the notation of an
asymptotic cone.

Asymptotic cones of groups capture global geometric properties of
Cayley graphs of these groups. For example, a hyperbolic group
``globally" looks like a tree because all asymptotic cones of them
are trees. In fact all asymptotic cones of all non-virtually cyclic
hyperbolic groups are isometric \cite{DP}.

Note that every asymptotic cone of a finitely generated group is a
complete homogeneous geodesic metric space. The asymptotic cones of
a group do not depend on the choice of observation points. So we
shall omit it from the notation. The subgroup $G^\omega_e$
consisting of all sequences $(g_i)$ in the ultrapower $G^\omega$
with $|g_i|=O(d_i)$ acts isometrically and transitively on the
asymptotic cone $\Con^\omega(G,(d_n))$ by left multiplication.

There are several basic questions about asymptotic cones that are
still unanswered (see below).

\subsection{The number of asymptotic cones}

Having unique asymptotic cone or at least asymptotic cones that do
not differ much is very convenient because in the applications (see
below), one never knows which asymptotic cone of a group one gets
because the scaling constants $d_n$ are almost never defined
explicitely. Unfortunately, the situation with the number of
distinct asymptotic cones of a group is far from ideal.

In \cite{KSTT}, it is proved (answering a question of Gromov) that
any uniform lattice $\Gamma$ in ${\mathrm SL}_n(\R)$, $n\ge 3$, has
$2^{2^{\aleph_0}}$ pairwise non-homeomorphic asymptotic cones
provided the Continuum Hypothesis (CH) is {\em not } true. If CH is
true, then $\Gamma$ has only one asymptotic cone. Moreover if CH is
true then any group has at most continuum asymptotic cones (roughly,
that is because by a theorem of Shelah, every asymptotic cone is
determined by its ``elementary theory", and there are at most
continuum different elementary theories). In \cite{DS}, we
constructed a group with continuum pair-wise non-$\pi_1$-equivalent
asymptotic cones (without any assumption about CH). The group is
given by an infinite presentation satisfying a small cancelation
condition (similar to the construction of ``monsters" from the
previous section). So it is not finitely presented, not amenable,
etc. The following problem can be probably attributed to Gromov.

\begin{problem} Is there a finitely presented group with continuum
asymptotic cones (assuming CH)?
\end{problem}

The answer is probably ``yes". The group from \cite{DS} is
recursively presented, so by Higman's theorem, it can be embedded
into a finitely presented group. There are many different versions
of Higman's embedding \cite{BORS}, \cite{OSamen}, \cite{OScol}, etc.
One can try to prove that there exists an embedding preserving the
property of having continuum non-$\pi_1$-equivalent asymptotic
cones.

Note that the first (independent of the Continuum Hypothesis)
example of a finitely presented group with two
non-$\pi_1$-equivalent cones is constructed in \cite{OS2005}.
Elementary amenable (infinitely presented) groups need not have
unique asymptotic cone \cite{OOS}.

Note that having different asymptotic cones means that
geometrically, a group looks different at different scales. Thus it
is possible that a ``less extreme" group always has unique
asymptotic cone (assuming CH). For example, Pansu \cite{Pansu}
proved that every nilpotent group has unique asymptotic cone.

\begin{problem} Is there a CAT(0)-group (a group acting properly co-compactly on a
finite dimensional CAT(0)-space) having several non-homeomorphic
asymptotic cones?
\end{problem}

%From \cite{Ham}, it seems to follow that every mapping class group
%of a surface has unique (up to bi-Lipschitz equivalence) asymptotic
%cone.

\begin{problem} Let $F$ be the R. Thompson group, $G$ be one of the Grigorchuk
groups of intermediate growth or any other self-similar group. How
many non-homeomorphic asymptotic cones does $F$ (resp. $G$) have?
\end{problem}

Paper \cite{OOS}  gives large classes of groups all of whose
asymptotic cones are locally isometric, but not all of them are
isometric. Similar methods can be used to show that the groups from
\cite{TV} and from \cite[Section 7]{DS} also satisfy this property.
However all these groups are infinitely presented. Moreover, in all
our examples, asymptotic cones are locally isometric to an $\mathbb
R$--tree, which implies hyperbolicity for finitely presented groups
\cite{OOS}. However the following problem is still open.

\begin{prob}[Olshanskii, Osin, Sapir]
Does there exist a finitely presented group all of whose asymptotic
cones are locally isometric, but not all of them are isometric?
\end{prob}

\begin{prob}[Olshanskii, Osin, Sapir]
Is there a non-virtually cyclic amenable group all of whose asymptotic
cones are locally isometric to a tree?
\end{prob}

An elementary amenable group with one asymptotic cone an $\R$-tree
was constructed in \cite{OOS} (that answered a question of B.
Kleiner).

\subsection{Fundamental groups of asymptotic cones}

Fundamental groups of asymptotic cones provide some important
information about the structure of the group. For example, if all
asymptotic cones are simply connected, then the group is finitely
presented, has polynomial Dehn function, and linear isodiametric
function (Gromov, \cite{Gr3}). Since any asymptotic cone is
homogeneous, uncountable, and its isometry group has usually
uncountable point stabilizers, every non-trivial loop in the
asymptotic cone has uncountably many copies having the same base
point. This was probably the motivation of Gromov's problem: Is
there an asymptotic cone of finitely generated group, whose
fundamental group is non-trivial and not of order continuum? Such a
group has been constructed in \cite{OOS}. In fact the asymptotic
cone of that group is homeomorphic to the direct product of an
$\R$-tree and a circle. So its fundamental group is $\Z$. Note that
the group in \cite{OOS} is not finitely presented, finding a
finitely presented example would be very interesting.

Since the asymptotic cones of a direct product are isometric to
direct products of asymptotic cones of the factors, one can realize
any finitely generated free Abelian group as the fundamental group
of $\CG $ for a suitable $G$ by taking direct products of groups. It
is quite possible that similarly one can construct an asymptotic
cone with a finite Abelian fundamental group.

\begin{problem}[Olshanskii, Osin, Sapir]
Does there exist a finitely generated group $G$ such that $\pi_1(\CG
)$ is countable (or, better, finitely generated) and non--Abelian
for some (any) $d$ and $\omega$? Can $\pi_1(\CG)$ be finite  and
non-Abelian?
\end{problem}

Note that for every countable group $C$ there exists a finitely
generated group $G$ and an asymptotic cone $\CG$ such that
$\pi_1(\CG)$ is isomorphic to the uncountable free power of $C$
\cite[Theorem 7.33]{DS}.

\section{Actions of groups on tree-graded spaces}
\label{tgs}

One of the most important applications of asymptotic cones is the
following observation due to Bestvina and Paulin: if a group
$\Gamma$ has infinitely many pairwise non-conjugate in a group $G$
homomorphisms $\phi\colon \Gamma\to G$, then $\Gamma$ acts
non-trivially by isometries on an asymptotic cone of $G$. The action
is the following: $$\gamma \circ (x_i)=(\phi_i(\gamma)x_i).$$ More
generally, if a group admits ``many" actions by isometries on a
metric space $X$, then it acts non-trivially on an asymptotic cone
of $X$. Thus, for example, if a hyperbolic group $G$ has infinite
$Out(G)$, then $G$ acts non-trivially on an $\R$-tree.

If $G$ is not hyperbolic, then asymptotic cones of $G$ need not be
trees. But they often are {\em tree-graded spaces}.

Recall \cite{DS} that a geodesic metric space $T$ is {\em
tree-graded with respect to a collection of geodesic subsets} $\pp$
(called {\em pieces} if any two subsets in $\pp$ intersect by at
most a point and every simple loop in $T$ is contained in one of the
pieces. General properties of tree-graded spaces can be found in
\cite{DS}.

In \cite{DS1}, it is proved that a group acting on a tree-graded
space non-trivially (i.e. not stabilizing a piece and not fixing a
point) also acts on an $\R$-tree without a global fixed points.
Moreover, the stabilizers of points and arcs of the $\R$-tree are
described.

Groups with tree-graded asymptotic cones occur much more often than
groups whose asymptotic cones are $\R$-trees. It is proved in
\cite{DS} that a group $G$ is relatively hyperbolic with respect to
proper subgroups $H_1,...,H_n$ if and only if is asymptotic cones
are tree-graded with respect to the collection of all
$\omega$-limits of sequence of cosets of the subgroups
$H_1,...,H_n$. Also it is noticed in \cite{DS}, that a geodesic
metric space is tree-graded with respect to a collection of proper
pieces if and only if it has a global cut-point. Such a situation
occurs very often. For example, if a group is given by a graded
small cancelation presentation \cite{book,Ols,OOS}, then its
asymptotic cones are {\em circle-trees}, that is tree-graded with
respect to embedded circles of diameters bounded from below and from
above \cite{OOS}. Asymptotic cones of mapping class groups are
tree-graded \cite{Behr}. So are the asymptotic cones of acylindrical
amalgamated products of groups \cite{DMS}. See more examples in
\cite{Behr,DS1,DMS}.

Using the result of Bestvina and Paulin, and the Rips theory of
groups acting on trees (developed further by Bestvina, Feighn,
Levitt, Sela, Guirardel and others), we found \cite{DS1} many
properties of subgroups of relatively hyperbolic groups that are
similar to properties of subgroups of hyperbolic groups. For
example, if a subgroup $H$ of a relatively hyperbolic group $G$ has
property (T),  then ${\mathrm Aut}(H)$ is commensurable with the
normalizer of $H$ in $G$. For other properties of subgroups of
relatively hyperbolic groups obtained by using actions on
tree-graded spaces see \cite{DS1}.

\begin{problem} Can one extend results of \cite{DS1} to subgroups of
other groups with tree-graded asymptotic cones?
\end{problem}

According to \cite{DS1}, in order to do that, one needs to:

\begin{itemize}
\item Describe the pieces of the tree-graded structure of the
asymptotic cones of $G$;
\item Describe the stabilizers of the action of the isometry group $G^\omega_e$
\begin{itemize}
\item a piece;
\item a pair of distinct pieces;
\item a pair (a piece and a point in it);
\item a pair of points not in the same piece.
\end{itemize}
\end{itemize}

 It certainly can be done for a much
larger class than the class of relatively hyperbolic groups. For
example, Olshanskii, Osin and I found results similar to the results
in \cite{DS1} for subgroups of groups given by {\em graded small
cancelation presentations} \cite{OOS}. That class includes torsion
groups, groups with all proper subgroups cyclic, and other monsters,
in particular.

It would be extremely interesting to extend results of \cite{DS1} to
subgroups of mapping class groups. The pieces of the asymptotic
cones of the mapping class groups have been described by Behrstock,
Kleiner, Minsky, and Mosher. The work by Behrstock, Dru\c tu and
myself on completing the program is currently in progress.

\section{Lacunary hyperbolic groups}

Lacunary hyperbolic groups were introduced in \cite{OOS}. A group is
called {\em lacunary hyperbolic} if one of its asymptotic cones is
an  $\R$-tree. By \cite{OOS}, a finitely generated group $G=\la
S\ra$ is lacunary hyperbolic if and only if it is a direct limit of
groups $G_i=\la S\ra$ and surjective homomorphisms $\alpha_i\colon
G_i\to G_{i+1}$ that is an identity on $S$, such that each $G_i$ is
$\delta_i$-hyperbolic, each $\alpha_i$ is injective on a ball of
radius $r_i$ and $\delta_i=o(r_i)$.

Every finitely presented lacunary hyperbolic group is hyperbolic
(Kapovich-Kleiner, \cite{OOS}). But that class is much larger than
the class of hyperbolic groups: Tarski monsters, torsion groups, and
other ``monsters" discussed above, can be (and many known examples
are by construction) lacunary hyperbolic.  There are amenable
non-virtually cyclic lacunary hyperbolic groups (I have already
mentioned that above) \cite{OOS}.

Although the class of lacunary hyperbolic groups is very large,
groups in that class share many common properties with hyperbolic
groups (see \cite{OOS}). For example,
\begin{itemize}
\item an undistorted subgroup of a \ch group is \ch itself,
\item a \ch group cannot contain a copy of $\Z^2$, an infinite
finitely generated subgroup of bounded torsion and exponential
growth, a copy of the lamplighter group, etc.,
\item  every \ch group is embedded into a relatively hyperbolic 2-generated
\ch group as a peripheral subgroup,
\item any
group that is hyperbolic relative to a \ch subgroup is \ch itself.
\end{itemize}

There are probably more such properties. Here are two concrete
problems.

\begin{prob}[Olshanskii,Osin,Sapir] Is it true that the growth of every non-elementary
\ch group is (a) exponential? (b) uniformly exponential?
\end{prob}

\begin{prob}[Olshanskii, Osin, Sapir]\label{expgr} Can a finitely generated non virtually
cyclic subgroup of exponential growth of a \ch group satisfy a
non-trivial law?
\end{prob}

\begin{prob} Is there an analog of Bestvina-Feighn combination
theorem \cite{BF} for \ch groups?
\end{prob}

In \cite{OOS}, it is proved that a subgroup with a non-trivial law
of a \ch group cannot have relative exponential growth.

The answer to Problem \ref{expgr} is ``no" for \ch groups for which,
using the notation of the definition of \ch groups, the injectivity
radii $r_i$ are ``much larger" than the hyperbolicity constants
$\delta_i$. More precisely it is enough to assume that
$r_i=\exp\exp(\delta_i)$ \cite{OOS}.

All known ``monsters" mentioned in Section \ref{s1} can be
constructed in such a way that this growth condition is satisfied.

It is also interesting to study linearity of \ch groups. We do not
know the answer to the following basic question.

\begin{prob}[Olshanskii, Osin, Sapir] Is every linear \ch group hyperbolic?
\end{prob}

\section{Linear and residually finite hyperbolic groups}

It is known (Kapovich, \cite{Kap1}) that a hyperbolic group can be
non-linear. But the existing examples are very implicit (although
their presentations can be, in principle, written down). I think
that non-linear hyperbolic groups may appear very often. One of the
simple sources of hyperbolic groups are ascending HNN extensions of
free groups. Let $F_n=\la x_1,...,x_n\ra$, $n\ge 2$, be the free
group, $\phi\colon F_n\to F_n$ be an injective homomorphism,
$\HNN_\phi(F_n)$ be the corresponding HNN extension.

It is known that $\HNN_\phi(F_n)$ is always coherent \cite{FH}
(every finitely generated subgroup is finitely presented) and
residually finite \cite{BS}. Some of these groups are not linear
(example: $\la x,y, t\mid txt\iv =x^2, tyt\iv = y^2\ra$
\cite{DS05}). But known examples all contain Baumslag-Solitar
subgroups, so these examples are not hyperbolic. Note that most
groups $\HNN_\phi(F_n)$ are small cancelation, hence hyperbolic.

\begin{problem}\label{45} Are there non-linear hyperbolic groups of the form
$\HNN_\phi(F_n)$ for non-surjective $\phi$? In particular, is the
group $\la x,y,t\mid txt\iv =xy, tyt\iv=yx\ra$ linear?
\end{problem}

It is known (Minasyan) that this group is hyperbolic. If $\phi$ is
surjective (i.e. an automorphism), then $\HNN_\phi(F_2)$ is always
linear. This follows from the linearity of the braid group $B_4$
(see \cite{Form}). It is quite possible (but still unknown) that
$\HNN_\phi(F_n)$ is always linear if $\phi$ is an automorphism.

The conjectural answers to the first question in Problem \ref{45} is
``yes" and to the second question ``no". Note that some groups of
the form $\HNN_\phi(F_n)$ are even inside $\SL_2(\R)$. The first
example is of course the Baumslag-Solitar group $BS(1,n)$ for every
$n$. But there are more complicated examples found in \cite{CD}. All
known examples are non-hyperbolic (contain Baumslag-Solitar
subgroups), have $n\ge 5$, and correspond to reducible endomorphisms
$\phi$. Are there hyperbolic ascending HNN extensions
$\HNN_\phi(F_2)$ with non-surjective but injective $\phi$ inside
$\SL_2(\R)$?

If one considers two injective but not surjective homomorphisms
$\phi, \psi\colon F_n\to F_n$ then one can form the {\em double} HNN
extension $\HNN_{\phi,\psi}(F_n)$ which is also quite often
hyperbolic. Methods from \cite{BS} (based on studying periodic
points of polynomial maps over finite fields) do not give a proof
that such groups are residually finite. To the contrary, they
indicate that these groups, generically, are {\em not} residually
finite because pairs of independent polynomial maps over a finite
field should not have common periodic points.

Thus it is possible that a double HNN extension of a free group can
be hyperbolic and non-residually finite.

\begin{problem} Are there non-residually finite hyperbolic groups of
the form $\HNN_{\phi,\psi}(F_n)$? In particular, is the group $$\la
x,y, t, u\mid txt\iv =xy, tyt\iv =yx, uxu\iv = [x,y], uyu\iv =
[x^2y,y^2x]\ra$$ residually finite (this is just a random choice of
``independent" endomorphisms $\phi, \psi$)?
\end{problem}

The conjectural answer is ``yes" to the first question and ``no" to
the second.

\section{R. Thompson group and other diagram groups}
\label{stg}

\subsection{Diagram groups}  One of the definitions of diagram groups is the following
(see \cite{GS}). Consider an alphabet $X$ and a set ${\cal S}$ of
{\em cells}, each cell is a polygon whose boundary is subdivided
into two directed paths (the top path and the bottom path having
common initial and terminal points) labeled by positive words $u$
(the top path) and $v$ (the bottom path) in the alphabet $X$. One
can consider the cell as a rewriting rule $u\to v$. Each cell $\pi$
is an {\em elementary $(u,v)$-diagram} with top path labeled by $u$,
bottom path labeled by $v$, and two distinguished vertices $\iota$
and $\tau$: the common starting and ending points of the top and
bottom paths. For every $x\in X$, there exists also the {\em
trivial} $(x,x)$-diagram: an edge labeled by $x$. Its top path and
bottom path coincide. There are four operations allowing to
construct more complicated diagrams from the elementary ones. These
are defined as follows.

%TeXCAD Options
%\grade{\on}
%\emlines{\off}
%\epic{\off}
%\beziermacro{\on}
%\reduce{\on}
%\snapping{\off}
%\quality{8.00}
%\graddiff{0.01}
%\snapasp{1}
%\zoom{4.0000}
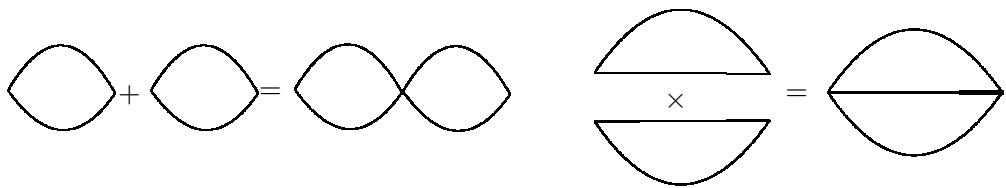
\unitlength .5mm % = 1.42pt
\linethickness{0.4pt}
\ifx\plotpoint\undefined\newsavebox{\plotpoint}\fi % GNUPLOT compatibility
\begin{picture}(269.75,64.88)(0,0)
\qbezier(5.25,26.5)(19.5,50.88)(33.75,25.75)
\qbezier(43.25,26.5)(57.5,50.88)(71.75,25.75)
\qbezier(81.5,26.75)(95.75,51.13)(110,26)
\qbezier(110.25,26.25)(124.5,50.63)(138.75,25.5)
\qbezier(33.75,25.75)(20,5.88)(5.25,26.5)
\qbezier(71.75,25.75)(58,5.88)(43.25,26.5)
\qbezier(110,26)(96.25,6.13)(81.5,26.75)
\qbezier(138.75,25.5)(125,5.63)(110.25,26.25)
\put(37.5,25.5){\makebox(0,0)[cc]{$+$}}
\put(75,26){\makebox(0,0)[cc]{$=$}}
\qbezier(161.25,31.25)(184,64.88)(207.75,31)
\qbezier(223.25,26)(246,59.63)(269.75,25.75)
\qbezier(161.25,18.25)(184,-15.37)(207.75,18.5)
\qbezier(223.25,26)(246,-7.62)(269.75,26.25)
%\emline(207.75,31)(161.5,31.25)
\multiput(207.75,31)(-11.5625,.0625){4}{\line(-1,0){11.5625}}
%\end
%\emline(269.75,25.75)(223.5,26)
\multiput(269.75,25.75)(-11.5625,.0625){4}{\line(-1,0){11.5625}}
%\end
%\emline(207.75,18.5)(161.5,18.25)
\multiput(207.75,18.5)(-11.5625,-.0625){4}{\line(-1,0){11.5625}}
%\end
%\emline(269.75,26.25)(223.5,26)
\multiput(269.75,26.25)(-11.5625,-.0625){4}{\line(-1,0){11.5625}}
%\end
%\emline(161.5,31.25)(161.5,31.25)
\put(161.5,31.25){\line(0,1){0}}
%\end
\put(182.75,23.75){\makebox(0,0)[cc]{$\times$}}
\put(215,25.5){\makebox(0,0)[cc]{$=$}}
\end{picture}

\begin{itemize}
\item The addition: $\Delta_1+\Delta_2$ is obtained by identifying the
terminal vertex $\tau$ of $\Delta_1$ with the initial vertex $\iota$
of $\Delta_1$. The top and the bottom paths of $\Delta_1+\Delta_2$
are defined in a natural way.

\item The multiplication: If the label of the bottom path of
$\Delta_1$ coincides with the label of the top path $\Delta_2$, then
$\Delta_1\Delta_2$ is defined by identifying the bottom path of
$\Delta_1$ with the top path of $\Delta_2$.

\item The inversion: $\Delta\iv$ is obtained from $\Delta$ by switching
the top and the bottom paths of the diagram.

\item Dipole cancelation: if $\pi$ is an $(u,v)$-cell, then we
identify $\pi\pi\iv$ with the trivial $(u,u)$-diagram. Thus we can
always replace a subdiagram $\pi\pi\iv$ by the trivial
$(u,u)$-subdiagram.
\end{itemize}

For every word $u$, the set of all $(u,u)$-diagrams forms a group
under the product operation, the {\em diagram group} with base word
$u$ and the given collection of cells. The $(u,u)$-diagrams are
called {\em spherical}.

\begin{example}[\cite{GS}] Here are some examples of diagram groups.

\begin{itemize}
\item The R. Thompson group $F$ is the diagram group of all
$(x,x)$-diagrams corresponding to the 1-letter alphabet $\{x\}$ and
one cell $x^2\to x$.

\item The wreath product $\Z\wr \Z$ is the diagram group of
$(ac,ac)$-digrams over the alphabet $$\{a, b_1, b_2, b_3, c\}$$
corresponding to cells $ab_1\to a, b_1\to b_2, b_2\to b_3. b_3\to
b_1, b_1c\to c$.

\item The free group $F_2$ is the diagram group of $(a, a)$-diagrams
over the alphabet $$\{a,a_1,a_2, a_3, a_4\}$$ and cells $$a\to a_1,
a_1\to a_2, a_2\to a, a\to a_3, a_3\to a_4, a_4\to a.$$

\item The direct product $\Z\times \Z$ is the diagram group of
$(ab,ab)$-diagram over the alphabet $$\{a, a_1, a_2, b, b_1, b_2\}$$
and cells $$a\to a_1, a_1\to a_2, a_2\to a, b\to b_1, b_1\to b_2,
b_2\to b.$$

\item Many right angled Artin groups are diagram groups \cite{GS3}.
\end{itemize}

\end{example}

The class of diagram groups is closed under direct and free products
\cite{GS1}, each diagram group is linearly orderable \cite{GS3} (and
so it is torsion-free). One can view a diagram group as a
2-dimensional analog of a free group (a free group is the group of
1-paths of a graph; the diagram groups are groups of 2-paths on
directed 2-complexes \cite{GS2}). The word problem in any subgroup
of a diagram group is very easy to decide. In many important cases
(including the Thompson group $F$), the conjugacy problem in a
diagram group has also an easy diagrammatic solution.

Although the class of diagram groups is large, some basic facts are
still missing.

\begin{problem}[Guba, Sapir]\label{hg} Is there a hyperbolic non-free diagram group?
\end{problem}

It is not even easy to construct a hyperbolic non-free subgroup of a
diagram group (i.e. representable by diagrams). Nevertheless it is
proved in \cite{CSS} that fundamental groups of the orientable
surface of genus 2 and the non-orientable surface of genus 4 are
representable by diagrams. The fundamental group of a non-orientable
surface cannot be a diagram group because all integral homology
groups of all diagram groups are free Abelian \cite{GS2}.

There is a reason to believe that the answer to problem \ref{hg} is
negative: by \cite{GS}, the centralizer of a diagram in a diagram
group is always a direct product of cyclic groups and diagram groups
(over the same collection of cells but with different base words).
The number of factors is equal to the maximal number of summands in
a representation of a conjugate of the diagram (by $(u,v)$-diagrams)
as a sum of spherical summands. Hence if a diagram group does not
contain $\Z\times \Z$, then all diagrams in that group and all their
conjugates are indecomposable into sums of spherical diagrams. So
far the only diagram groups satisfying this property were free.

\subsection{The R. Thompson group $F$ and its subgroups}

The R. Thompson group $F$, corresponding to the ``simplest"
collection of cells $\{xx\to x\}$, plays an important role in the
class of diagram groups. It is proved in \cite{GS3} that every
diagram group is a subgroup of the so-called universal diagram group
$U$. The group $U$ is finitely presented, and is explicitly
constructed as a split extension of a right angled Artin group and
the group $F$.

Subgroups of R.Thompson's group $F$ have been of great interest
after the first result of Thompson and Brin-Squier \cite{BrSq} that
$F$ contains no free non-Abelian subgroups. There are several
Tits-like dichotomies concerning subgroups of $F$:

\begin{itemize}
\item (Brin-Squier) A subgroup of $F$ is either Abelian or contains a copy of
$\Z^\infty$.
\item (Guba-Sapir) A subgroup of $F$ is either Abelian or contains a
copy of $\Z\wr\Z$.
\item (Bleak) A subgroup of $F$ is either solvable or contains a
copy of the (restricted) direct product $W$ of wreath products $\Z
\wr \Z \wr...\wr Z$ ($n$ times), $n=1,2,...$.
\end{itemize}

Bleak's papers \cite{Bleak1,Bleak2} contain, in a sense, a complete
description of all solvable subgroups of $F$. Note that nilpotent
subgroups of diagram groups are free Abelian \cite{GS}.

The following problem of Brin (also suggested by Guba and me) is
still open and very interesting.

\begin{problem}[Brin]\label{pbrin} Is it true that every subgroup of $F$ either
contains a copy of $F$ or is elementary amenable.
\end{problem}

It is easy to see that $F$ is not elementary amenable itself (since
$F'$ is simple). Note that Brin's examples \cite{Brin1} show that
even finitely generated elementary amenable subgroups subgroups of
$F$ can be very complicated. The most promising methods related to
Problem \ref{pbrin} seems to be those developed by Brin
\cite{Brin1,Brin2} and Bleak \cite{Bleak1,Bleak2}. But the diagram
group methods may prove fruitful too.

One of the strongest results about the group $F$ is Guba's theorem
that $F$ has quadratic Dehn function \cite{Guba}. (Hence by
Papasoglu's result \cite{Pap}, the asymptotic cones of $F$ are
simply connected.) The following question  remains open and even
more intriguing after \cite{Guba}.

\begin{problem}[Guba, Sapir] Is $F$ automatic?
\end{problem}

\subsection{Amenability, property A, and embeddings into Hilbert spaces}

Of course, I should also mention one of the central problems about
R.Thompson's group:

\begin{problem}[R. Thompson, R. Geoghegan] Is $F$ amenable?
\end{problem}

The closest property to amenability that $F$ provably has is
a-T-menability: existence of a proper action on a Hilbert space
(Farley, \cite{Farley}). It is not known even if $F$ satisfies the
much weaker than amenability G. Yu's property A.

A metric space $X$ satisfies property A if for every $\epsilon>0$,
and every $R$, there exists a collection of finite subsets
$A_x\subset X\times \N$, $x\in X$, and a number $S>0$ such that

\begin{itemize}
\item For every $x,y\in X$ with $\dist(x,y)<R$, we have
$$\frac{|A_x\Delta A_y|}{|A_x\cap A_y|}<\epsilon,$$ where $\Delta$
denotes the symmetric difference,
\item For every $x\in X$, the diameter of the projection of $A_x$
onto $X$ is at most $S$.
\end{itemize}
 G. Yu proved that if a Cayley graph $X$ of a group $G$
satisfies property $A$, then it coarsely (uniformly) embeds into a
Hilbert space $H$, i.e. there exists a function $f\colon X\to H$ and
an increasing unbounded {\em compression} function $\rho\colon
\R_+\to\R_+$ such that
$$\rho(\dist(x,y)) \le  \dist(f(x), f(y))\le \dist(x,y)$$
for every $x,y\in X$.

Guentner and Kaminker \cite{GK} proved that, ``conversely", if such
an embedding exists and, in addition $$\lim_{n\to\infty}
\frac{\rho(n)}{\sqrt{n}}=\infty$$ (i.e. $\sqrt{n}\ll\rho(n)$), then
$X$ satisfies property $A$. If, moreover, the embedding $f$ is {\em
equivariant}, that is for some proper action of $G$ on $H$, we have
$f(gx)=gf(x)$ (for all $g,x\in G$), then $G$ is amenable. In
\cite{AGS}, we constructed (using actions of diagram groups on the
so-called 2-trees, and a result of Burillo \cite{Bur1}) an
equivariant embedding of $F$ into a Hilbert space with
$$\liminf_{n\to \infty} \frac{\rho(n)}{\sqrt{n}}>0$$ and
proved that for every (not necessary equivariant) coarse embedding
of $F$ into a Hilbert space, $$\limsup_{n\to\infty}
\frac{\rho(n)}{\sqrt{n}\log n} <\infty.$$ Thus there is a gap
$[\sqrt{n}, \sqrt{n}\log n]$ where a possible compression function
of $F$ may belong. Hence there is still hope that the answer to the
following problem is positive and $F$ has property $A$ or even is
amenable.

\begin{problem}[Arzhantseva, Guba,Sapir]\label{ags} Is there
(equivariant or not) coarse embedding of the R.Thompson group $F$
into a Hilbert space with compression function $\rho$ satisfying
$\lim_{n\to\infty} \frac{\rho(n)}{\sqrt{n}}=\infty$?
\end{problem}

In \cite{ADS}, we introduced the notion of a {\em Hilbert space
compression} gap of a metric space $X$.

\begin{definition}\label{defcomprfct}
Let $(X,\dist )$ be a metric space. Let $f, g\colon\R_+ \to \R_+$ be
two increasing functions such that $f\ll g$ and $\lim_{x\to \infty }
f(x)= \lim_{x\to\infty} g(x)=\infty .$

We say that $(f,g)$ is a  {\em Hilbert space compression gap of
$(X,\dist )$ } if

\begin{itemize}
\item[(1)]
 there exists an embedding of $X$ into
a Hilbert space with compression function at least $f$;

\item[(2)]
for every embedding of $X$ into a Hilbert space, its compression
function $\rho$ satisfies $\rho\ll g$.

\end{itemize}
\end{definition}

The quotient $g/f$ is called the {\em size} of the compression gap.
We showed in \cite{ADS}, in particular, that for every real $\alpha$
between $0$ and $1$, there exists a finitely generated group of
asymptotic dimension at most $2$ such that for every $\epsilon>0$,
$\left(\frac{x^\alpha}{\log^{1+\epsilon}(x+1)}\, ,\,
x^\alpha\right)$ is a Hilbert space compression gap of that group
(that answered a question of Guentner and Niblo). Results of
\cite{ADS} show that Thompson's group $F$ has compression gap of
logarithmic size. Similar results for free groups, lattices in Lie
groups, etc. follow from results of Bourgain \cite{Bourgain} and
Tessera \cite{Tessera}.

\begin{problem} Is there any group with compression gap less than
$\log x$ (say, $\log\log x$)?
\end{problem}

The answer would be positive (for the R. Thompson group) if Problem
\ref{ags} had positive answer.

Note that every finitely generated diagram group such as $F$, $\Z\wr
\Z$, etc. has two natural metrics: the word metric and the diagram
metric. The distance between two diagrams $\Delta_1$ and $\Delta_2$
in the diagram metric is the number of cells in the diagram
$\Delta_1\iv \Delta_2$ after removing all the dipoles. Translated
into the language of diagram groups, the result of Burillo
\cite{Bur1} used in \cite{AGS} is the following: for the group $F$,
the word metric and the diagram metric are quasi-isometric. In
\cite{AGS}, we proved that the same is true for $\Z\wr\Z$ and one of
the universal diagram groups.

\begin{problem}[Arzhantseva, Guba, Sapir]
Does every finitely generated diagram group satisfy
the Burillo property?
\end{problem}

If the answer is ``yes", then every finitely generated diagram group
would have an embedding into a Hilbert space with compression
function $\rho(n)\succ \sqrt{n}$.

\subsection{The membership problem and distortion of subgroups}
\label{dist}

Many algorithmic problems of $F$ are known to be  solvable. Those
are the word problem, the conjugacy problem \cite{GS}, the multiple
conjugacy problem \cite{KM} and others. One of the most important
problems that are still open is the {\em generalized word problem},
that is the membership problem in the finitely generated subgroups.

\begin{definition}\label{ddist} Recall that an increasing function $f(n)\colon \N\to\N$ is
called the distortion function of a finitely generated subgroup
$H=\la Y\ra$ of a group $G=\la X\ra$ if every element $h\in H$ that
has length $n$ in $G$ (relative to $X$) has length at most $f(n)$ in
$H$ (relative to $Y$), and $f$ is the smallest function with this
property.
\end{definition}

If we change generating sets $X$ and $Y$ (keeping them finite), the
new distortion function is {\em equivalent} to the old one. We say
that increasing functions $f$ and $g$ are equivalent if
$$\frac1Cf(\frac nC)-C <g(n)<Cf(Cn)+C$$ for some $C\ge 1$ and all
sufficiently large $n$.

It is well known (probably Farb was the first who noticed that) that
the generalized word problem is solvable for a given finitely
generated subgroup $H$ of a group $G$ if and only if the distortion
of $H$ is recursive. Very little is known about the distortion of
subgroups of $F$. It is known that $F$ and its ``brothers" $F_n$,
direct powers $F^n$ for every $n$, $F\wr \Z$, etc. (see
\cite{GS,GS1}) can be embedded into $F$ without distortion. It is
much more complicated to find distorted subgroups of $F$. There are
known examples of subgroups with non-linear distortion \cite{GS1}.
But all these subgroups are solvable and the distortion is at most
quadratic.

\begin{problem} [Guba-Sapir]\label{Brin6} a) Does $F$ contain a distorted copy of
itself (note that by Brin \cite{Brin2}, $F$ contains {\em lots} of
copies of itself).

b) Is there a finitely generated subgroup of $F$ with bigger than
polynomial (or even bigger than quadratic) distortion?

c) Is there a finitely generated subgroup of $F$ with non-recursive
distortion?
\end{problem}

Problem \ref{Brin6}a) has been also suggested by Brin.

\section{Surface subgroups of right angled Artin groups}
\label{ssg}

\subsection{The main problem}

This is another kind of representation problems. We want to
represent closed hyperbolic surface groups (and more generally
hyperbolic non-free groups) as subgroups of right angled Artin
groups. Right angled Artin groups have so many nice properties
(surveyed in Charney \cite{Ch}) that embedding a group into one of
them is certainly interesting. There is an industry developed mostly
by D. Wise and F. Haglund  designed to embed a group into {\em a}
right angled Artin group. Here we are interested in a different kind
of problems: what subgroups a given right angled Artin group can
have?

\begin{problem}[Crisp-Sageev-Sapir] Is there an algorithm to decide given a finite graph
$K$ whether the right angled Artin group $A(K)$ contains the
fundamental group of a closed hyperbolic surface?
\end{problem}

There are several partial results related to this problem
\cite{DSS}, \cite{Kim}, \cite{CSS}. Clearly if $K_1$ is a full
subgraph of $K_2$, then $A(K_1)\le A(K_2)$, so if $A(K_1)$ contains
$\pi_1(S)$ for some surface $S$ then so does $A(K_2)$. Kim
\cite{Kim} proved that if $K_1$ is obtained from $K$ by the {\em
co-contraction} of a pair of non-adjacent vertices, then $A(K_1)$ is
inside $A(K)$. The {\em co-contraction} of a pair $(a,b)$ amounts to
replacing $a, b$ by a new vertex $c$ connected to all the common
neighbors of $a, b$ in $K$.

Hence it is enough to describe graphs $P$ such that no proper
subgraph $P'$ of $P$ and no graph $P'$ obtained from $P$ by
co-contraction is such that $A(P')$ contains a hyperbolic closed
surface subgroup. Let us call such graphs {\em extremal}. There is
one series of extremal graphs: the $n$-cycles for $n\ge 5$
\cite{DSS}. There are also eight exceptional extremal graphs
$P_1(6), P_2(6), P_1(7), P_2(7), P_1(8)-P_4(8)$ (the number in
parentheses is the number of vertices in the graph)  found in
\cite{CSS} (Kim also found $P_1(6)$, the triangular prism). Is it
the full collection of extremal graphs? We know \cite{CSS} that this
collection contains all extremal graphs with up to $8$ vertices.

There are also several reduction moves $P\to P'$ found in \cite{CSS}
such that if $A(P)$ contains a hyperbolic closed surface subgroup,
then so does $P'$ (and $P'$ contains fewer vertices or edges than
$P$). We do not know if this collection of reduction moves is
complete, namely we do not know whether the following statement is
true: If $A(P)$ does not contain a hyperbolic closed surface
subgroup, then $A(P')$ does not as well, for one of our reduction
moves $P\to P'$. This collection of moves is complete for 8-vertex
graphs \cite{CSS}.

\subsection{Dissection diagrams}

Let $S$ be a surface (possibly with boundary). Let $G=\la X\mid
R\ra$ be a finitely presented group. Let $\Psi$ be a van Kampen
diagram over the presentation of $G$ drawn on $S$. Roughly speaking
it is a graph drawn on $S$ with edges labeled by letters from $X$,
such that each connected component of $S\setminus \Psi^1$ is a
polygon with boundary path labeled by a word from $R^{\pm 1}$ (see
more details in \cite{book}).

Given a van Kampen diagram $\Psi$ on $S$, one can define a
homomorphism $\psi\colon \pi_1(S)\to G$ as follows. As a base-point,
pick a vertex $v$ of $\Psi$. Let $\gamma$ be any loop at $v$. Since
all cells in the tessellation $\Psi$ are polygons, $\gamma$ is
homotopic to a curve that is a composition of edges of $\Psi$. Then
$\psi(\gamma)$ is the word obtained by reading the labels of edges
of $\Psi$ along $\gamma$. Since the label of the boundary of every
cell in $\Psi$ is equal to 1 in $G$, the words corresponding to any
two homotopic loops $\gamma$, $\gamma'$ represent the same element
in $G$. Hence $\psi$ is indeed well-defined. The fact that $\psi$ is
a homomorphism is obvious.

Conversely, the standard argument involving $K(.,1)$-complexes gives
that every injective homomorphism $\psi\colon \pi_1(S)\to G$
corresponds in the above sense to a van Kampen diagram over $G$ on
$S$.

If $G=A(K)$ is a right angled Artin group, then every cell in a van
Kampen diagram is a square, and instead of a van Kampen diagram on
$S$, it is convenient to consider its dual picture (pick a point
inside every cell, connect the points in neighbor cells by an edge
labeled by the label of the common edge of the cells). It is called
the {\em $K$-dissection diagram} of the surface, and was introduced
by Crisp and Wiest in \cite{CW}. The edges having the same labels
form collections of pairwise disjoint simple closed orientation
preserving curves and arcs connecting points on the boundary of $S$.
Each of these curves has a natural transversal direction. Each curve
is labeled by a vertex of $K$, two curves intersect only if their
labels are adjacent in $K$.

If $\Delta$ is the $K$-dissection diagram corresponding to a van
Kampen diagram $\Psi$ on $S$, then the corresponding homomorphism
$\psi\colon \pi_1(S)\to A(K)$ takes any loop $\gamma$ based at $v$
to the word of labels of the dissection curves and arcs of $\Delta$
crossed by $\gamma$ (a letter in the word can occur with exponent
$1$ or $-1$ according to the direction of the dissection curve
crossed by $\gamma$). There are several partial algorithms allowing
to check whether a homomorphism $\psi$ corresponding to the
$K$-dissection diagram is injective \cite{CSS}. But the answer to
the next question is still unknown.

\begin{problem}[Crisp-Sageev-Sapir] Is there an algorithm which
given a $K$-dissection diagram on a surface $S$, decides whether the
corresponding homomorphism $\psi\colon\pi_1(S)\to A(K)$ is
injective?
\end{problem}

\section{Surface subgroups of hyperbolic groups and ascending HNN
extensions of free groups} \label{sse}

A very interesting problem by Gromov (on the Bestvina problem list)
asks whether every non-virtually free hyperbolic group contains a
hyperbolic closed surface subgroups. Here is a concrete partial case
of this problem.

\begin{problem} Does the group $\la x,y,t \mid txt\iv =xy,
tyt\iv=yx\ra$ contain a hyperbolic closed surface subgroup?
\end{problem}

Of course a negative answer is more interesting than a positive one.

According to the previous section, to solve this problem, one would
need to study van Kampen diagrams over this presentation on
surfaces. Here is a more general algorithmic problem.

\begin{problem} For which injective but not surjective
endomorphisms $\phi\colon F_k\to F_k$, the HNN extension
$\HNN_\phi(F_k)$ contains a hyperbolic closed surface subgroup?
\end{problem}

Note that subgroups of such HNN extensions have been described in
details by Feighn and Handel in \cite{FH}.

\section{Percolation in Cayley graphs of groups\\ by Iva Koz\'akov\'a and Mark Sapir}
\label{sp}

\subsection{Basic properties of percolation}

Geometry of random structures associated with Cayley graphs of
groups is very interesting in general. There is a lot of work
related, say, to random walks on Cayley graphs (starting with the
classical cases of Cayley graphs of $\Z$ and $\Z\times \Z$ with the
natural generators), see Woess, \cite{Woess}. Here we shall
formulate some problems related to random subgraphs of Cayley graphs
and percolation. For more problems, motivation and discussion,  see
the survey of Benjamini and Schramm \cite{BeSch}.

\begin{definition}
A {\em Bernoulli bond percolation} on $\mathcal{G}$ is a product
probability measure $\mathrm{P}_p$ on the space $\Omega=\{0,1\}^E$,
the subsets of the edge set $E$. For $0\leq p\leq 1$ the product
measure is defined via $\mathrm{P}_p(\omega(e)=1)=p$ for all $e\in
E$.
\end{definition}

An element $\omega$ of the probability space $\Omega$ is called {\em
configuration} or {\em realization} of percolation. For any  $\omega
\in \Omega$, the bond $e\in E$ is called {\em open} if $\omega(e)=1$
and {\em closed} otherwise. Thus each bond is open with probability
$p$ independently of all other bonds.

We write $\mathrm{E}_p$ for the {\em expected value} with respect to
$\mathrm{P}_p$. For any configuration $\omega$, open edges form a
random subgraph of $\mathcal{G}$.

\begin{definition}
An {\em(open) cluster} is a connected component of such subgraph
$\omega$. An open cluster containing the origin is denoted by $C$
and the number of vertices in $C$ by $|C|$.  The {\em percolation
function} is defined to be the probability that the origin is
contained in an infinite cluster, i.e.
$\theta(p)=\mathrm{P}_p(|C|=\infty)$.
\end{definition}

The connectivity function $\tau_p(o,x)=\mathrm{P}_p(o\leftrightarrow
x)$ is the probability that there is an open path connecting
vertices $o$ and $x$. We use $\chi(p)=\mathrm{E}_p(|C|)$ for the
mean size of the open cluster at the origin.

As  $p$ grows from $0$ to $1$, the behavior of the percolation
process changes. We can distinguish several phases of similar
characteristic.

\begin{definition} For any graph, we define three {\em critical probabilities}:
\begin{align*}
 p_c&=\sup\{p: \theta(p)=0\},\\
p_u&=\inf\{p:\text{ There is a unique infinite cluster }\mathrm{P}_p\text{-almost surely}\}, \\
p_\textrm{exp}&=\sup\{p : \exists_{C,\gamma>0} \forall_{x,y \in V}
\tau_p(x,y)\leq C e^{-\gamma \mathrm{dist}(x,y)} \}.
\end{align*}
\end{definition}

H\"aggstr\"om, Peres and Schramm \cite{HPS} showed that for any
homogeneous graph, for every $p$ such that $0\leq p < p_c$ all
clusters are finite, there are infinitely many infinite clusters if
$p_c< p < p_u$, and if $p_u < p \leq1$, then there is unique
infinite cluster $\mathrm{P}_p$-almost surely.

It is known \cite{Sch} that $p_c\le p_{exp}\le p_u$ for every Cayley
graph. It is interesting to characterize all groups (resp. Cayley
graphs) such that $p_c$ and $p_u$ are different. It is known that
$p_c=p_u$ for all amenable groups. The converse has not been proved
yet.

\begin{problem}[Benjamini-Schramm \cite{BeSch}] Is it true that $p_c\ne p_u$ for
every Cayley graph of a non-amenable group?
\end{problem}

The conjectural answer is ``yes" \cite{BeSch}. There exist partial
results in this direction. Pak and Smirnova-Nagnibeda \cite{Pak}]
proved that for any nonamenable (finitely generated) group there
exists a finite symmetric set of generators $S$ in $G$ such that for
the Cayley graph of $G$ with respect to $S$, $p_c<p_u$. Schonmann
\cite{Sch} showed that $p_c<p_u$ for highly nonamenable groups (i.e.
such that the edge-isoperimetric (Cheeger) constant is bigger than
$(\sqrt{2d^2-1}-1)/2$, where $d$ is the vertex degree). Gaboriau
\cite{Gab} proved that the first $\ell^2$-Betti number of a group
does not exceed $\frac{1}{2} (p_u-p_c)$ for every Cayley graph of
$G$, therefore as soon as $\ell^2$-Betti number is nonzero,
$p_u>p_c$.

The study of percolation started with the square lattice $\Z^2$.
Using the self-duality of this lattice it can be shown that
$p_c=\frac{1}{2}$ \cite{Ke}, although the proof is very non-trivial.
The exact value of $p_c$ for $d$-dimensional cubic lattices $\Z^d$,
$d\ge 3$, is unknown, only numerical approximations are available.

A well understood situation is the case of a regular tree, i.e. the
Cayley graph of a free group with respect to the free generating
set. For a regular tree, the critical probability
$p_c=\frac{1}{d-1}$, where $d$ is the vertex degree.

Using the tree-graded structure of Cayley graphs of free products
Koz\'akov\'a \cite{Ko} proved that for a free product of groups
$G_1*G_2$ with the natural generating set, the $p_c$ is a root of
$(\chi_1(p)-1)(\chi_2(p)-1)=1$, where $\chi_i(p)$ is the expected
cluster size in the Cayley graph of $G_i$, \cite{Ko}. For example
for ${\rm PSL}_2(\Z)=\Z/2\Z * \Z/3\Z$, $p_c$ is $.5199...$, the
unique root of the polynomial $2p^5-6p^4+2p^3+4p^2-1$ in the
interval $(0,1)$.

A degenerated case is the Cayley graph of the infinite cyclic group
$\Z$, or any finite extension of it. For these graphs, $p_c=1$. And
again the converse is not known in general. The question can be
formulated in the following way using the growth function.

\begin{problem}[Benjamini, Schramm \cite{BeSch}] Assume that
$G$ has growth function faster than linear, is $p_c<1$?
\end{problem}

The conjectural answer is ``yes" \cite{BeSch}. Only amenable case is
of interest, since for every nonamenable graph, the
edge-isoperimetric (Cheeger) constant $i_E$ is positive and
$p_c<\frac{1}{i_E+1}$ \cite{BeSch}. It has been shown by Lyons (see
\cite{BeSch2}) that for groups with polynomial or exponential growth
$p_c<1$. The same is true for all finitely presented groups
\cite{BaBe}. All Grigorchuk groups of intermediate growth have
subgroups that are direct products of two proper infinite subgroups,
therefore $p_c<1$ (noticed by Muchnik and Pak \cite{MP}).

It is also not known whether the properties $p_c=1$ or $p_c=p_u$ are
invariant under the change of generators of a group or even the same
for any graph that is quasi-isometric to the Cayley graph of a group
(physicists believe that it is the case at least for $\Z^d$).

\begin{problem}[Smirnova-Nagnibeda] What are $p_c, p_u, p_{exp}$ for $\SL_2(\Z)$, the
R.Thompson group $F$, Grigorchuk groups (and other automata groups)?
\end{problem}

It may be difficult to find precise values, but even observations
whether $p_c, p_u, p_{exp}$ are all different are interesting.

On planar graphs a very useful property is that any bi-infinite
simple path splits the graph into two components (for examples of
applications of this see \cite{BeSch3} and \cite{BaBe}). Can a
similar property hold for other graphs?

\begin{problem} Is it true that for a hyperbolic group an infinite
cluster ($p>p_c$) separates the graph $P_p$-a.s.?
\end{problem}

An open cluster as a subgraph of the Cayley graph carries some
properties of the original graph, but its geometry can be also
different. Gaboriau \cite{Gab}  related  harmonic Dirichlet
functions on a graph to those on the infinite clusters in the
uniqueness phase. One way to describe the change in the geometry of
a subgraph is the distortion (compare with Definition \ref{ddist}).

\begin{definition}
The {\em distortion} function of a connected subgraph $C$
(containing origin $O$) of a graph $\mathcal{G}$ is given by
\begin{align*}
D(n) &= \max_{\stackrel{y\in C}{\dist_{\mathcal{G}}(O,y)\leq n}}
\dist_C(O,y)
\end{align*}
where $\dist_C$ (resp. $\dist_\mathcal{G}$) is the graph metric in
$C$ (resp. in $\mathcal{G}$). In the context of percolation, $D(n)$
gives the distortion of an open cluster $C$ (at the origin) for any
given realization $\omega\in\Omega$, and we define the {\em expected
distortion} as $\mathrm{E}_p(D(n))$.
\end{definition}

Note that different definitions of distortion of clusters related to
random walks in random environments are also useful (see \cite{BLS},
for example).

For $p<p_c$ the distortion is bounded by the cluster size, therefore
the expected distortion is a constant (asymptotically). On trees the
distance in a subgraph does not differ from the original one,
therefore if $p>p_c$ the distortion of an infinite open cluster on a
tree is linear. The result by Antal and Pisztora \cite{AP} suggests
that the expected distortion for the $d$-dimensional cubic lattice
$\Z^d$ should be at most linear as well, they proved that there
exists $\rho=\rho(d,p)$ such that $\mathrm{P}_p$-almost surely
$$
\limsup_{\dist_{\Z^d}(O,y)\rightarrow
\infty}\frac{\dist_C(O,y)}{\dist_{\Z^d}(O,y)}<\rho.
$$

\begin{problem} What is the possible expected distortion of an open cluster
in a Cayley graph of a group? Same question for hyperbolic groups is
also open.
\end{problem}

The most interesting percolation phenomena occur, when the parameter
$p$ is near its critical value $p_c$.

The principal hypothesis of the percolation theory
%of the scaling theory
is that various quantities (like the probability of the cluster at
the origin being infinite, or the mean size of this cluster) have
certain specific asymptotic behavior near the critical point $p_c$.

\begin{definition}\label{defex}
Assume $\theta(p)$ is continuous at $p_c$ and that
\begin{align*}
\theta(p)&\approx(p-p_c)^\beta\quad \text{ as }p \searrow p_c,\\
\chi(p)&\approx(p_c-p)^{-\gamma}\quad \text{ as } p \nearrow p_c,\\
%\xi(p)&\approx(p_c-p)^{-\nu}\quad \text{ as }p \nearrow p_c,\\
\mathrm{P}_{p_c}(|C|=n)&\approx n^{-1-\frac{1}{\delta}}.
\end{align*}
Then we say that $\beta$, $\gamma$ and $\delta$ %and $\nu$
are {\em critical exponents} of the graph.
\end{definition}

For a regular tree it is possible to compute these exponents
directly. The critical exponents have the following so called {\em
mean-field values}:
\begin{align*}
\beta=1, \gamma=1, \delta=2. %, \nu=\frac{1}{2}.
\end{align*}

Physicists believe that the numerical values of critical exponents
depend only on the underlying space and not on the structure of the
particular lattice.

It was proved by Hara and Slade \cite{HaSl} that the critical
exponents of a $d$-dimensional cubic lattice take their mean-field
values for $d\geq 19$. It is believed to hold even for $d>6$. In
fact it might hold for all Cayley graphs of groups with growth
function $\ge n^6$ (compare with Problem \ref{ph} below).

%***copied from paper:*
A well known conjecture in percolation theory claims that the
critical exponents of all  lattices in ${\mathbb R}^2$ are the same
($\beta=5/36$, $\gamma=-43/18$) as proved for triangular lattices by
Smirnov and Werner \cite{SmWe}). The motivation for this conjecture
is that for every lattice $L$ in ${\mathbb R}^2$ the
Gromov-Hausdorff limit of rescaled copies of $L$, $L/2$, $L/3$,
\dots is isometric to ${\mathbb R}^2$ with the $L_1$-metric, i.e.
any two lattices in ${\mathbb R}^2$ have the same asymptotic cones
(\cite{Gr3}). Sapir conjectured that this is true in general: if two
groups have isometric asymptotic cones corresponding to the same
ultrafilters and the same scaling constants, then their critical
exponents should coincide. In particular, since all asymptotic cones
of all non-elementary Gromov-hyperbolic groups are isometric (they
are isometric to the universal ${\mathbb R}$-tree of degree
continuum by a result of \cite{DP}), every non-elementary hyperbolic
group should have mean-field valued critical exponents. Thus we
formulate

\begin{problem}\label{ph} Is it true that the critical exponents of all Cayley
graphs of groups with isometric asymptotic cones are equal? In
particular, is it true that every Cayley graph of a non-elementary
hyperbolic group has mean-field valued critical exponents?
\end{problem}

By \cite{Sch}, the critical exponents take their mean-field values
for all non-amenable planar graphs with one end, and for unimodular
graphs with infinitely many ends (in particular, for all Cayley
graphs of groups with infinitely many ends).

\subsection{Scale invariant groups}

One of the methods used by physicists in studying percolation in
$\Z^2$ (and $\Z^d$ in general) is {\em rescaling} or {\em
renormalization}. For some $n>1$, we consider a tessellation of
$\Z^d$ by cubes of size $n$. That tessellation can be viewed as a
new (rescaled) copy $\Gamma=n\Z$ of $\Z^d$. Given a percolation
realization on $\Z^d$, we say that an edge $e$ of the $n\times
n$-cube in $\Gamma$ is {\em open} if one can cross the cube in the
direction of the edge using open edges of $\Z^d$ (i.e. there is an
open path connecting the $(d-1)$-faces of the cube that are
orthogonal to $e$).

This collection of open and closed edges in $\Gamma$ can be viewed
(under certain assumptions that are obvious to physicists) as a
realization of percolation on $\Gamma$ (see Grimmett \cite{Grim}).
If the original realization had an infinite cluster, the new one
would have infinite clusters too. This way, starting with a
super-critical percolation ($p>p_c$) we get a new supercritical
percolation with a bigger $p$. Hence one is able to deduce
information about $p_c$ and critical exponents \cite{Grim}.

The rescaling method is based on the fact that $\Z^d$ contains lots
of copies of itself of finite index. In general, following
Benjamini, we call a finitely generated group $G$ {\em scale
invariant} if it has subgroups of finite index that are isomorphic
to $G$ and the intersection of all such subgroups of finite index is
finite.

\begin{problem} [Benjamini,
http://www.math.weizmann.ac.il/$\sim$itai/] Is it true that every
scale invariant finitely generated group is virtually nilpotent?
\end{problem}

One way to deal with this problem would be to consider the actions
of the group on itself induced by the finite index embeddings of the
group into itself. That gives an action of the group on its
asymptotic cone (as in \cite{DS1} and Section \ref{tgs}). This
should imply, in particular, that a relatively hyperbolic group is
scale invariant only if it is elementary.

Note that R. Thompson group $F$ has many copies of itself as
subgroups of finite index. But $F$ is not residually finite, hence
not scale invariant. Grigorchuk's groups have many subgroups of
finite index which are isomorphic to finite direct powers of the
ambient group. But these subgroups are not isomorphic to the ambient
group, so Grigorchuk groups are not scale invariant too.

{\bf Acknowledgement.} I am very grateful to Iva K\'ozak\'ova for
explaining percolation theory to me and for co-writing Section 9 of
the paper. I am also grateful to I. Benjamini, M. Brin, R. Lyons, T.
Smirnova-Nagnibeda and A. Olshanskii for their comments and
suggestions.

\begin{minipage}[t]{2.6 in}
\noindent Mark V. Sapir\\ Department of Mathematics\\
Vanderbilt University\\
m.sapir@vanderbilt.edu\\
%http://www.math.vanderbilt.edu/$\sim$msapir\\
\end{minipage}

\end{document}